\documentstyle[11pt]{article}
 \newcommand{\qb}[2]{{\left [{#1 \atop #2} \right]}}
 \newlength{\standardunitlength}
\newtheorem{cor}{Corollary} \newtheorem{lemma}{Lemma}
\newtheorem{theorem}{Theorem} \newtheorem{prop}{Proposition}
\newenvironment{proof}{\noindent {\sc Proof:}}{$\Box$ \vspace{2 ex}}

\begin{document}

\begin{center} Cellini's descent algebra, dynamical systems, and
semisimple conjugacy classes of finite groups of Lie type \end{center}

\begin{center}
By Jason Fulman
\end{center}

\begin{center}
Stanford University
\end{center}

\begin{center}
Department of Mathematics
\end{center}

\begin{center}
Building 380, MC 2125
\end{center}

\begin{center}
Stanford, CA 94305
\end{center}

\begin{center}
email:fulman@math.stanford.edu
\end{center}

1991 AMS Primary Subject Classifications: 20G40, 20F55

\newpage
Proposed running head: Cellini's descent algebra

\newpage \begin{abstract} By algebraic group theory, there is a map
from the semisimple conjugacy classes of a finite group of Lie type to
the conjugacy classes of the Weyl group. Picking a semisimple class
uniformly at random yields a probability measure on conjugacy classes
of the Weyl group. We conjecture that this measure agrees with a
second measure on conjugacy classes of the Weyl group induced by a
construction of Cellini which uses the affine Weyl group. This is
verified in some cases such as type $C$ odd characteristic. For the
identity conjugacy class in type $A$, the proof of the conjecture
amounts to an interesting number theoretic reciprocity law. More
generally the type $A$ case leads to number theory involving Ramanujan
sums. Models of card shuffling, old and new, arise naturally. In type
$C$ even characteristic connections are given with dynamical
systems. We indicate, at least in type $A$, how to associate to a
semisimple conjugacy class an element of the Weyl group, refining the
map to conjugacy classes. \end{abstract}


\section{Introduction}

	In performing a definitive analysis of the
Gilbert-Shannon-Reeds model of card-shuffling, Bayer and Diaconis
\cite{BD} defined a one-parameter family of probability measures on
the symmetric group $S_n$ called $k$-shuffles. Given a deck of $n$
cards, one cuts it into $k$ piles with probability of pile sizes
$j_1,\cdots,j_k$ given by $\frac{{n \choose
j_1,\cdots,j_k}}{k^n}$. Then cards are dropped from the packets with
probability proportional to the pile size at a given time (thus if the
current pile sizes are $A_1,\cdots,A_k$, the next card is dropped from
pile $i$ with probability $\frac{A_i}{A_1+\cdots+A_k}$). They proved
that $\frac{3}{2}log_2(n)$ $2$-shuffles are necessary and suffice to
mix up a deck of $n$ cards (\cite{F6} shows that the use of cuts does
not help to speed things up). These $k$-shuffles induce a probability
measure on conjugacy classes of $S_n$, hence on partitions $\lambda$
of $n$. Diaconis, McGrath, and Pitman \cite{DMP} studied the
factorization of random degree $n$ polynomials over a field $F_q$ into
irreducibles. The degrees of the irreducible factors of a randomly
chosen degree $n$ polynomial also give a random partition of $n$. The
main result of \cite{DMP} is that this measure on partitions of $n$
agrees with the measure induced by card shuffling when $k=q$. The
cycle structure of biased riffle shuffles was studied in \cite{F}.

	It is worth observing that the GSR measures are well-studied
and appear in many mathematical settings. The paper \cite{H} is a good
reference to applications to Hochschild homology (tracing back to
Gerstenhaber and Schack \cite{Ger}), and the paper
\cite{BergeronWolfgang} describes the relation with explicit versions
of the Poincar\'e-Birkhoff-Witt theorem. Section 3.8 of \cite{SSt}
describes the type $A$ measure in the language of Hopf algebras. In
recent work, Stanley \cite{Sta} has related biased riffle shuffles
with representation theory of the symmetric group, thereby giving an
elementary probabilistic interpretation of Schur functions and a
different approach to some work in the random matrix community. He
recasts and extends many of the results of \cite{BD} and \cite{F}
using quasisymmetric functions.

	The paper \cite{F2}, building on the papers \cite{BBHT} and
\cite{BHR}, defined card shuffling measures $H_{W,x}$ for any finite
Coxeter group and real non-zero $x$ (and actually for any real
hyperplane arrangement). These measures extend the GSR shuffles and
have attractive properties, the most important of which is that
$H_{W,x}(w) \geq 0$ if $W$ is crystallographic and $x$ is a good prime
for $W$. Observe that monic degree $n$ polynomials over $F_q$ are
precisely the semisimple conjugacy classes of $GL(n,q)$, and that
there is for general types a natural map $\Phi$ from semisimple
conjugacy classes of finite groups of Lie type to conjugacy classes of
the Weyl group (this will be reviewed in Section \ref{background}). It
would be natural to conjecture that the measures $H_{W,x}$ are related
to semisimple conjugacy classes of finite groups of Lie type. This is
not the case, and it emerges \cite{F3} that the measures $H_{W,x}$ are
instead related to the semisimple adjoint orbits of finite groups of
Lie type on their Lie algebras. By a theorem of Steinberg \cite{St}
the number of semisimple adjoint orbits and number of semisimple
conjugacy classes are both equal to $q^r$ where $r$ is the rank of the
group; however combinatorially they are quite different. One of the
aims of this paper is to understand this difference, giving analogs of
\cite{F3} for semisimple conjugacy classes.

	It will be supposed for the rest of the paper that $G$ is a
connected, semisimple, simply connected group defined over a finite
field of $q$ elements. Let $\cal G$ be the Lie algebra of $G$. Let $F$
denote both a Frobenius automorphism of $G$ and the corresponding
Frobenius automorphism of $\cal G$. Suppose that $G$ is
$F$-split. Thus in type $A$ the group under consideration is $SL(n,q)$
rather than $GL(n,q)$. The theory of how random characteristic
polynomials of $SL(n,q)$ factor into irreducibles is much more
complicated than the $GL(n,q)$ case, and involves for instance the
value of $n$ mod $q$.

	Let $\Pi=\{\alpha_1,\cdots,\alpha_r\}$ be a set of simple
roots for a root system of a finite Coxeter group $W$ of rank
$r$. Recall that the descent set $Des(w) \subseteq \Pi$ of an element
$w$ of $W$ is the set of simple positive roots which $w$ maps to
negative roots. Let $d(w)=|Des(w)|$. Solomon \cite{Sol} constructed
his so called descent algebra as follows. For any $I \subseteq \Pi$,
let \[ Y_I = \{w \in W| Des(w) = I\}. \] He proved that any product
$Y_J Y_K$ is a linear combination of the $Y_I$'s with non-negative
integral coefficients. Cellini \cite{Ce1} defined a different type of
descent algebra, replacing the set of simple roots $\Pi$ by all
vertices of the extended Dynkin diagram (thus including the
possibility that the highest root may be sent to a negative
root). This definition will be recalled in Section \ref{background},
and Section \ref{mod} indicates how in types $A$ and $C$ it gives rise
to some new physical models of card shuffling.

	Conjecture 1 relating semisimple conjugacy classes to
Cellini's work will be given in Section \ref{main}, together with its
proof in some cases. Explicit formulas for $x_k$ will be given for
type $A$. Even in the simplest possible case--the identity conjugacy
class in type $A$--Conjecture 1 is interesting and amounts to a number
theoretic reciprocity law. Section \ref{dynamical} describes
connections with dynamical systems. As one example, it is shown that a
proof of Conjecture 1 for type $C$ in even characteristic would give
an alternate solution to a problem concerning the cycle structure of
unimodal maps. Along the way. the Rogers-Weiss \cite{RW} enumeration
of transitive unimodal maps is given a first combintarial proof. In
Section \ref{refine}, we discuss the possibility of refining the map
$\Phi$ so as to associate to a semisimple conjugacy class an {\it
element} of the Weyl group. As with Cellini's work, the ideas rely
heavily on the affine Weyl group. We illustrate the discussion for the
symmetric group on three symbols and pose a more general conjecture.

	To close the introduction, we mention the follow-up paper
\cite{F6} which compares the type $A$ affine shuffles considered here
with ordinary riffle shuffles followed by cuts. Using representation
theoretic work on the Whitehouse module, the cycle strucutre of riffle
shuffles followed by a cut was determined. When $gcd(n,q-1)=1$, strong
evidence is given that these measures, though different, coincide at
the level of conjugacy classes.

\section{Background} \label{background}

	Recall that $G$ is a connected, semisimple, simply connected
group defined over a finite field of $q$ elements. Letting $F$ be a
Frobenius automorphism of $G$, we suppose that $G$ is $F$-split. Next
we recall the map $\Phi$ from semisimple conjugacy classes $c$ of
$G^F$ to conjugacy classes of the Weyl group. Since the derived group
of $G$ is simply connected (the derived group of a simply connected
group is itself), Theorem 3.5.6 of \cite{C} gives that the
centralizers of semisimple elements of $G$ are connected. Consequently
$x$ is determined up to conjugacy in $G^F$ and $C_G(x)$, the
centralizer in $G$ of $x$, is determined up to $G^F$ conjugacy. Let
$T$ be a maximally split maximal torus in $C_G(x)$. Then $T$ is an
$F$-stable maximal torus of $G$, determined up to $G^F$ conjugacy. By
Proposition 3.3.3 of \cite{C}, the $G^F$ conjugacy classes of
$F$-stable maximal tori of $G$ are in bijection with conjugacy classes
of $W$. Define $\Phi(c)$ to be the corresponding conjugacy class of
$W$.

	For example (e.g. page 273 of \cite{Mac}) in type $A_{n-1}$
the semisimple conjugacy classes of $SL(n,q)$ correspond to monic
degree $n$ polynomials $f(z)$ with constant term 1. Such a polynomial
factors as $\prod_i f_i^{a_i}$ where the $f_i$ are irreducible over
$F_q$. Letting $d_i$ be the degree of $f_i$, $\Phi(c)$ is the
conjugacy class of $S_n$ corresponding to the partition $(d_i^{a_i})$.

	Next we recall the work of Cellini \cite{C} (the definition
which follows differs slightly from hers, being inverse and also
making use of her Corollary 2.1). We follow her in supposing that $W$
is a Weyl group (i.e. crystallographic). Letting $\alpha_0$ denote the
negative of the highest root, let $\tilde{\Pi}=\Pi \cup
\alpha_0$. Define the cyclic descent $Cdes(w)$ to be the elements of
$\tilde{\Pi}$ mapped to negative roots by $w$, and let
$cd(w)=|Cdes(w)|$. For instance for $S_n$ the simple roots with
respect to a basis $e_1,\cdots,e_n$ are $e_i-e_{i+1}$ for
$i=1,\cdots,n-1$ and $\alpha_0=e_n-e_1$. Thus the permutation $4 \ 1 \
3 \ 2 \ 5$ (in 2-line form) has 3 cyclic descents.

	For $I \subseteq \tilde{\Pi}$, put \[ U_I = \{w \in W|Cdes(w)
\cap I = \emptyset\} .\] Let $Y$ be the coroot lattice. Then define
$a_{k,I}$ by

\[ \left\{ \begin{array}{ll}
 |\{ t \in Y| <\alpha_0,t> = k, <\alpha_i,t> = 0 \ for \ \alpha_i
\in I-\alpha_0, <\alpha_i,t> >0 \ for \ \alpha_i \in \tilde{\Pi}-I\}| &
\mbox{if $\alpha_0 \in I$}\\
 |\{ t \in Y| <\alpha_0,t> < k, <\alpha_i,t> = 0 \ for \ \alpha_i \in
I, <\alpha_i,t> >0 \ for \ \alpha_i \in \Pi-I\}| &
\mbox{if $\alpha_0 \not \in I$} \end{array}\right.\] Finally, define an
element $x_k$ of the group algebra of $W$ by \[ x_k = \frac{1}{k^r}
\sum_{I \subseteq \tilde{\Pi}} a_{k,I} \sum_{w \in U_I} w .\]
Equivalently, the coefficient of an element $w$ in $x_k$ is \[
\frac{1}{k^r} \sum_{I \subseteq \tilde{\Pi} - Cdes(w)} a_{k,I}.\]

	For the purpose of clarity, in type $A$ this says that the
coefficient of $w$ is $x_k$ is equal to $\frac{1}{k^{n-1}}$ multiplied
by the number of integers vectors $(v_1,\cdots,v_n)$ satisfying the
conditions

\begin{enumerate}
\item $v_1+\cdots+v_n=0$
\item $v_1 \geq v_2 \geq \cdots \geq v_n, v_1-v_n \leq k$
\item $v_i>v_{i+1}$ if $w(i)>w(i+1)$ (with $1 \leq i \leq n-1$)
\item $v_1<v_n+k$ if $w(n)>w(1)$
\end{enumerate}

	From Cellini (loc. cit.), it follows that the $x_k$ satisfy
the following two desirable properties:

\begin{enumerate}
\item (Measure) The sum of the coefficients in the expansion of $x_k$
in the basis of group elements is 1. Equivalently,

\[ \sum_{I \subseteq \tilde{\Pi}} a_{k,I} |U_I| = k^r .\] In
probabilistic terms, the element $x_k$ defines a probability measure
on the group $W$.

\item (Convolution) $x_k x_h = x_{kh}$.
\end{enumerate}

	The above definition of $x_k$ is computationally convenient
for this paper. We note that Cellini (loc. cit.) constructed the $x_k$
in the following more conceptual way, when $k$ is a positive
integer. Let $W_k$ be the index $k^r$ subgroup of the affine Weyl
group generated by reflections in the hyperplanes corresponding to
$\{\alpha_1,\cdots,\alpha_r\}$ and also the hyperplane
$\{<x,\alpha_0>=k\}$. There are $k^r$ unique minimal length coset
representatives for $W_k$ in the affine Weyl group, and $x_k$ is
obtained by projecting them to the Weyl group.

	The following problem is very natural. It would also be
interesting (at least in type $A$) to relate the $x_k$'s to cyclic and
Hochschild homology.

{\bf Problem:} Determine the eigenvalues (and multiplicities) of $x_k$
acting on the group algebra by left multiplication. More generally,
recall that the Fourier transform of a probability measure $P$ at an
irreducible representation $\rho$ is defined as $\sum_{w \in W} P(w)
\rho(w)$. For each $\rho$, what are the eigenvalues of this matrix?

\section{Physical Models of Card Shuffling} \label{mod}

	We pause to give examples which both illustrate the definition
of $x_k$ and gives a relation with models of card shuffling (some of
them new). Writing $x_k=\sum c_w w$ in the group algebra, the notation
$x_k^{-1}$ will denote $\sum c_w w^{-1}$.

\begin{prop} When $W$ is the symmetric group $S_{2n}$, the element
$(x_2)^{-1}$ has the following probabilistic interpretation:

	Step 1: Choose an even number between $1$ and $2n$ with the
probability of getting $2j$ equal to $\frac{{2n \choose
2j}}{2^{2n-1}}$. From the stack of $2n$ cards, form a second pile of
size $2j$ by removing the top $j$ cards of the stack, and then putting
the bottom $j$ cards of the first stack on top of them.

	Step 2: Now one has a stack of size $2n-2j$ and a stack of
size $2j$. Drop cards repeatedly according to the rule that if stacks
$1,2$ have sizes $A,B$ at some time, then the next card comes from
stack $1$ with probability $\frac{A}{A+B}$ and from stack 2 with
probability $\frac{B}{A+B}$. (This is equivalent to choosing uniformly
at random one of the ${2n \choose 2j}$ interleavings preserving the
relative orders of the cards in each stack).

	 The description of $x_2^{-1}$ is the same for the symmetric
group $S_{2n+1}$, except that at the beginning of Step 1, the chance
of getting $2j$ is $\frac{{2n+1 \choose 2j}}{2^{2n}}$ and at the
beginning of Step 2, one has a stack of size $2n+1-2j$ and a stack of
size $2j$. \end{prop}

\begin{proof} We argue for the case $S_{2n}$, the case of $S_{2n+1}$
being similar. Recall that in type $A_{2n-1}$ the coroot lattice is
all vectors with integer components and zero sum with respect to a
basis $e_1,\cdots,e_{2n}$, that $\alpha_i=e_i-e_{i+1}$ for
$i=1,\cdots,2n-1$ and that $\alpha_0=e_{2n}-e_1$. The elements of the
coroot lattice contributing to some $a_{2,I}$ are:

\[ \begin{array}{ll} (0,0,\cdots,0,0) & I=\tilde{\Pi}-\alpha_0\\
(1,0,0,\cdots,0,0,-1) & I=\tilde{\Pi}-\{\alpha_1,\alpha_{2n-1}\}\\
(1,1,0,0,\cdots,0,0,-1,-1) &
I=\tilde{\Pi}-\{\alpha_2,\alpha_{2n-2}\}\\ \cdots & \cdots\\
(1,1,\cdots,1,0,0,-1,\cdots,-1,-1) &
I=\tilde{\Pi}-\{\alpha_{n-1},\alpha_{n+1}\}\\
(1,1,\cdots,1,1,-1,-1,\cdots,-1,-1) & I=\tilde{\Pi}-\alpha_n
\end{array} \] One observes that the inverses of the permutations in
the above card shuffling description for a given $j$ contribute to
$u_I$ where

\[ I =             \left\{ \begin{array}{ll}
\tilde{\Pi}-\alpha_0 & \mbox{if $2j=0$}\\				
\tilde{\Pi}-\{\alpha_k,\alpha_{2n-k}\} & \mbox{if $2j=2 min(k,2n-k)$}\\
\tilde{\Pi}-\alpha_n & \mbox{if $2j=2n$}\end{array}
			\right.  \] The total number of such
permutations for a fixed value of $j$ is ${2n \choose 2j}$, the number
of interleavings of $2n-2j$ cards with $2j$ cards preserving the
relative orders in each pile. Since $\sum_{j=0}^n {2n \choose
2j}=2^{2n-1}$, and $x_2$ is a sum of $2^{2n-1}$ group elements, the
proof is complete. \end{proof}

	The following problem is very natural.

{\bf Problem:} Find a physical description the elements $x_k$ in type
$A$ for integer $k>2$.

	As a second example, we describe the elements $x_k$ in type
$C_n$ (this was essentially done for $k=2$ in \cite{Ce2}). Lemma
\ref{com}, which follows easily from Theorem 1 of \cite{Ce2}, gives a
formula for $x_k$. Lemma \ref{com} also implies that for $k$ odd the
$x_{k}$ measure is equal to the measure $H_{C_n,k}$ of \cite{F2}.

\begin{lemma} \label{com} Let $d(w)$ and $cd(w)$ denote the number of
descents and cyclic descents of $w \in C_n$. Then the coefficient of
$w$ in $x_k$ is

\[ \begin{array}{ll}
\frac{1}{k^n} {\frac{k-1}{2}+n-d(w) \choose n} & \mbox{$k$ odd}\\
\frac{1}{k^n} {\frac{k}{2}+n-cd(w) \choose n} & \mbox{$k$ even}
\end{array}\] \end{lemma}

\begin{proof} For the first assertion, from Theorem 1 of \cite{Ce2},
the coefficient of $w$ in $x_k$ is \begin{eqnarray*} \frac{1}{k^n}
\sum_{l=d(w)}^n {\frac{k-1}{2} \choose l} {n-d(w) \choose l-d(w)} & =
& \frac{1}{k^n} \sum_{l=d(w)}^n {\frac{k-1}{2} \choose l} {n-d(w)
\choose n-l}\\ & = & \frac{1}{k^n} \sum_{l=0}^n {\frac{k-1}{2} \choose
l} {n-d(w) \choose n-l}\\ & = & \frac{1}{k^n} {\frac{k-1}{2}+n-d(w)
\choose n}.  \end{eqnarray*} The second assertion is similar and
involves two cases. \end{proof}

	Proposition \ref{rel} shows that the elements $x_k$ in type
$C$ arise from physical models of card-shuffling. The models which
follow were made explicit for $k=2$ in \cite{BD} and for $k=3$ in
\cite{BB}. The implied formulas for card shuffling resulting from
combining Propositions \ref{com} and \ref{rel} may be of interest.

\begin{prop} \label{rel} The element $x_k^{-1}$ in type $C_n$ has the
following description:

	Step 1: Start with a deck of $n$ cards face down. Choose $q$
numbers $j_1,\cdots,j_k$ multinomially with the probability of getting
$j_1,\cdots,j_k$ equal to $\frac{{n \choose
j_1,\cdots,j_k}}{k^n}$. Make $k$ stacks of cards of sizes
$j_1,\cdots,j_k$ respectively. If $k$ is odd, then flip over the even
numbered stacks. If $k$ is even, the flip over the odd numbered stacks.

	Step 2: Drop cards from packets with probability proportional
to packet size at a given time. Equivalently, choose uniformly at
random one of the ${n \choose j_1,\cdots,j_k}$ interleavings of the
packets.
\end{prop}

\begin{proof} The proof proceeds in several cases, the goal being to
show that the inverse of the above processes generate $w$ with the
probabilities in Lemma \ref{com}. We give details for one subcase--the
others being similar--namely even $k$ when $w$ satisfies
$cd(w)=d(w)$. (The other case for $k$ even is $cd(w)=d(w)+1$). The
inverse of the probabilistic description in the theorem is as follows:

	Step 1: Start with an ordered deck of $n$ cards face
down. Successively and independently, cards are turned face up and
dealt into one of $k$ uniformly chosen random piles. The even piles
are then flipped over (so that the cards in these piles are face
down).

	Step 2: Collect the piles from pile 1 to pile $k$, so that
pile 1 is on top and pile $k$ is on the bottom.

Consider for instance the permutation $w$ given in 2-line form by $-2\
3\ 1\ 4\ -6\ -5\ 7$. Note that this satisfies $cd(w)=d(w)$ because the
top card has a negative value (i.e. is turned face up). It is
necessary to count the number of ways that $w$ could have arisen from
the inverse description. This one does using a bar and stars argument
as in \cite{BD}. Here the stars represent the $n$ cards, and the bars
represent the $k-1$ breaks between the different piles. It is easy to
see that each descent in $w$ forces the position of two bars, except
for the first descent which only forces one bar. Then the remaining
$(k-1)-(2d(w)-1)=k-2d(w)$ bars must be placed among the $n$ cards as
$\frac{k-2d(w)}{2}$ consecutive pairs (since the piles alternate
face-up, face-down). This can be done in ${\frac{k}{2}+n-cd(w) \choose
n}$ ways, proving the result. \end{proof}

	Proposition \ref{rel} leads to a direct proof of the
convolution property in type $C$, arguing along the lines of the proof
of the convolution property for the type $A$ riffle shuffles of
\cite{BD} (which are different from the type $A$ shuffles considered
here).

\section{First Main Conjecture} \label{main}

	This section presents the main conjecture relating the
elements $x_k$ to semisimple conjugacy classes of finite groups of Lie
type, followed by evidence in its favor.

{\bf Conjecture 1:} Let $G$ be a connected, semisimple, simply
connected group defined over a finite field of $q$ elements. Letting
$F$ be a Frobenius automorphism of $G$, suppose that $G$ is
$F$-split. Let $c$ be a semisimple conjugacy class of $G^F$ chosen
uniformly at random. Then for all conjugacy classes $C$ of the Weyl
group $W$, \[ \sum_{w \in C} Probability(\Phi(c) = C) = \sum_{w \in C}
Coef. \ of \ w \ in \ x_q. \]

	To begin, we derive an expression for $x_k$ in type
$A_{n-1}$. For this recall that the major index of $w$ is defined by
$maj(w) = \sum_{i: 1 \leq i \leq n-1 \atop w(i)>w(i+1)} i$. The
notation $\qb{n}{k}$ denotes the $q$-binomial coefficient
$\frac{(1-q)\cdots(1-q^n)}{(1-q)\cdots (1-q^k) (1-q) \cdots
(1-q^{n-k})}$. Let $C_m(n)$ denote the Ramanujan sum $\sum_k
e^{\frac{2\pi ikn}{m}}$ where $k$ runs over all integers less than and
prime to $m$. The following lemma of Von Sterneck (see \cite{Ram} for
a proof in English) will be helpful.

\begin{lemma} \label{VS} (\cite{V}) The number of ways of expressing
$n$ as the sum mod $m$ of $k \geq 1$ integers of the set
$0,1,2,\cdots,m-1$ repetitions being allowed is \[ \frac{1}{m}
\sum_{d|m,k} {\frac{m+k-d}{d} \choose \frac{k}{d}} C_d(n).\]
\end{lemma}

\begin{theorem} \label{typeAform} In type $A_{n-1}$, the coefficient
of $w$ in $x_k$ is equal to any of the following:

\begin{enumerate}

\item $\frac{1}{k^{n-1}}$ multiplied by the number of partitions with
$\leq n-1$ parts of size at most $k-cd(w)$ such that the total number
being partitioned has size congruent to $-maj(w) \ mod \ n$.

\item $\frac{1}{k^{n-1}}$ multiplied by the number of partitions with
$\leq k-cd(w)$ parts of size at most $n-1$ such that the total number
being partitioned has size congruent to $-maj(w) \ mod \ n$.

\item \[ \begin{array}{ll} \frac{1}{n k^{n-1}} \sum_{d|n,k-cd(w)} {\frac{n+k-cd(w)-d}{d} \choose
\frac{k-cd(w)}{d}} C_d(-maj(w)) & \mbox{if \ $k-cd(w)>0$}\\
\frac{1}{k^{n-1}}  & \mbox{if $k-cd(w)=0, maj(w)=0$ mod n}\\
0 & \mbox{otherwise}
\end{array} \]

\item \[ \frac{1}{k^{n-1}} \sum_{r=0}^{\infty} Coeff. \ of \ q^{r \cdot n} \
in \ \left(q^{maj(w)} \qb{k+n-cd(w)-1}{n-1}\right).\] \end{enumerate} \end{theorem}

\begin{proof} Using the definition of $x_k$, one sees that 

\begin{eqnarray*}
 x_k & = & \frac{1}{k^{n-1}} \sum_{I \subseteq
\tilde{\Pi}-Cdes(w)} a_{k,I}\\
 & = & \sum_{v_1+\cdots+v_n=0, v_1 \geq \cdots v_{n-1} \geq v_n, v_1-v_n \leq k, \vec{v}
\in Z^n, \atop v_i
> v_{i+1} \ if \ e_i-e_{i+1} \in Cdes(w), \ and \ 
v_1<v_n+k \ if \ \alpha_0 \in Cdes(w)} 1\\
 & = & \frac{1}{k^{n-1}}
\sum_{r=0}^{\infty} Coeff. \ of \ q^{r \cdot n} \ in \sum_{ k \geq v_1
\geq \cdots v_{n-1} \geq v_n = 0, \vec{v} \in Z^n \atop v_i > v_{i+1}
\ if \ e_i-e_{i+1} \in Cdes(w), \ and \ v_1<v_n+k \ if \ 
\alpha_0 \in Cdes(w)} q^{\sum v_i}.
\end{eqnarray*} Now let $v_i'=v_i - |\{j: i \leq j \leq n-1, w(j)>w(j+1)\}|$. Then the expression for $x_k$ simplifies to

\begin{eqnarray*} 
&& \frac{1}{k^{n-1}}
\sum_{r=0}^{\infty} Coeff. \ of \ q^{r \cdot n} \ in \sum_{ k-cd(w) \geq v_1'
\geq \cdots v_{n-1}' \geq v_n' = 0, \vec{v} \in Z^n} q^{\sum v_i' + \sum_i |\{j: i \leq j \leq n-1, w(j)>w(j+1)\}|}\\
& = &  \frac{1}{k^{n-1}}
\sum_{r=0}^{\infty} Coeff. \ of \ q^{r \cdot n} \ in \ q^{maj(w)} \sum_{ k-cd(w) \geq v_1'
\geq \cdots v_{n-1}' \geq v_n' = 0, \vec{v} \in Z^n} q^{\sum v_i'}.
\end{eqnarray*}

	This proves the first assertion of the theorem. The second
assertion follows from the first by viewing partitions diagramatically
and taking transposes. The third assertion follows from the second and
Lemma \ref{VS}. The fourth assertion follows from either the first or
second assertions together with the well-known fact that the
generating function for partitions with at most $a$ parts of size at
most $b$ is the q-binomial coefficient $\qb{a+b}{a}$.\end{proof}

	Note that the second formula simplifies when $n=k$ is and $n$
is prime, since $1 \leq cd(w) <n$ for all $w \in S_n$. The follow-up
\cite{F6} proves Conjecture 1 in this case.

	Next we prove Conjecture 1 for the symmetric group $S_3$,
leaving the case $S_2$ as an exercise to the reader. Recall that the
semisimple conjugacy classes of $SL(3,q)$ are monic degree 3
polynomials with constant term $1$ and are $q^2$ in number.

\begin{theorem} \label{checkS3} Conjecture 1 holds for $S_3$.
\end{theorem}

\begin{proof} Theorem \ref{typeAform} calculates $x_q$ for $S_3$. Thus
it is only necessary to calculate the left hand side of the quantity
in Conjecture 1 for each conjugacy class of $S_3$. For the identity
conjugacy class this amounts to counting polynomials of the form
$(x-a)(x-b)(x-c)$ with $abc=1$. By the principle of inclusion and
exclusion, the number of such polyomials with distinct roots is
$\frac{(q-1)^2-3(q-1)+2e}{6}$ where $e$ is the number of cube roots of
$1$ in $F_q$. The number of solutions with exactly two of $\{a,b,c\}$
equal is $q-1-e$ and the number of solutions with $a=b=c$ is $e$. Thus
the total number of such polynomials is \[ \begin{array}{ll}
\frac{q^2+q}{6} & \mbox{if $q=0,2 \ mod \ 3$}\\ \frac{q^2+q+4}{6} &
\mbox{if $q=1 \ mod \ 3$} \end{array} \] agreeing with the coefficient
of the identity in $x_q$. For the conjugacy class corresponding to a
tranposition, one must count monic degree 3 polynomials factoring into
a quadratic and a linear factor, with product of the roots equal to
1. This is simply the number of degree 2 irreducible polynomials,
i.e. $\frac{q^2-q}{2}$, since the degree 1 factor is then
determined. This answer agrees with the answer in the first
step. Since the only class left is the class of 3-cycles and because
the total number of semisimple conjugacy classes of $SL(3,q)$ is
$q^2$, the proof is complete.  \end{proof}

	 As a consequence of Lemma \ref{VS} and Theorem
\ref{typeAform}, one obtains a verification of Conjecture 1 for the
identity conjugacy class in type $A$. For this we need the following
corollary of Lemma \ref{VS}. It can be regarded as a sort of modular
combinatorial reciprocity law. As such, a direct combinatorial proof
is desirable.

\begin{cor} \label{recip} For any positive integers $x,y$, the
number of ways (disregarding order and allowing repetition) of writing
$n$ (mod $y$) as the sum of $x$ integers of the set $0,1,\cdots,y-1$
is equal to the number of ways (disregarding order and allowing
repetition) of writing $n$ (mod $x$) as the sum of $y$ integers of the
set $0,1,\cdots,x-1$. \end{cor}

\begin{proof} It is enough to show that

\[ \frac{1}{y} \sum_{d|x,y} {\frac{x+y-d}{d} \choose \frac{x}{d}}
C_d(n) = \frac{1}{x} \sum_{d|x,y} {\frac{x+y-d}{d} \choose
\frac{y}{d}} C_d(n).\] Elementary algebra shows that the terms
corresponding to a given value of $d$ on both sides of this equation
are equal. \end{proof}

\begin{cor} \label{checkident} Conjecture 1 holds for the identity
conjugacy class in type $A$. \end{cor}

\begin{proof} The semisimple conjugacy classes $c$ of $SL(n,q)$ such
that $\Phi(c)=id$ are simply the number of ways (disregarding order)
of picking $n$ elements of $F_q^*$ whose product is the
identity. After choosing a generator for $F_q^*$ viewed as a cyclic
group of order $q-1$, one sees that the number of such semisimple
conjugacy classes is the number of ways of expressing $0$ as the sum
mod $q-1$ of $n$ integers from the set $0,1,2,\cdots,q-2$. By Lemma
\ref{VS} and the fact that there are $q^{n-1}$ semisimple conjugacy
classes, the probability that $\Phi(c)=id$ is \[
\frac{1}{(q-1)q^{n-1}} \sum_{d|q-1,n} {\frac{n+q-1-d}{d} \choose
\frac{n}{d}}.\] By Theorem \ref{typeAform} the coefficient of the identity in
$x_q$ is equal to \[\frac{1}{nq^{n-1}} \sum_{d|n,q-1}
{\frac{n+q-1-d}{d} \choose \frac{q-1}{d}}.\] The result follows by
arguing as in Corollary \ref{recip}. \end{proof}

	Proposition \ref{reformulation} shows that Conjecture 1 has an
reformulation in terms of generating functions for type $A$. This
reformulation is interesting because one side is mod $n$ and the other
side is mod $k-1$! For its proof, Lemma \ref{likeVic} will be helpful.

\begin{lemma} \label{likeVic} Let $f_{n,k,d}$ be the coefficient of
$z^n$ in $(\frac{z^k-1}{z-1})^d$. Then $1/i \sum_{d|i} \mu(d)
f_{m,k,i/d}$ is the number of size $i$ aperiodic necklaces on the
symbols $\{0,1,\cdots,k-1\}$ with total symbol sum $m$.  \end{lemma}

\begin{proof} This is an elementary Mobius inversion running along the
lines of a result in \cite{R}. \end{proof}

\begin{theorem} \label{reformulation} Let $f_{n,k,d}$ be the
coefficient of $z^n$ in $(\frac{z^k-1}{z-1})^d$. Let $n_i(w)$ be the
number of $i$-cycles in a permutation $w$. Then Conjecture 1 in type
$A$ is equivalent to the assertion (which we intentionally do not
simplify) that for all $n,k$,

\begin{eqnarray*}
& & \sum_{m=0 \ mod \ n} Coef. \ of \ q^m u^n t^k \
in \ \sum_{n=0}^{\infty} \frac{u^n} {(1-tq)\cdots(1-tq^n)} \sum_{w
\in S_n} t^{cd(w)} q^{maj(w)} \prod x_i^{n_i(w)}\\ & = &
\sum_{m=0 \ mod \ k-1} Coef. \ of \ q^m u^n t^k \ in \
\sum_{k=0}^{\infty} t^k \prod_{i=1}^{\infty} \prod_{m=1}^{\infty}
(\frac{1}{1-q^mx_iu^i})^{1/i \sum_{d|i} \mu(d) f_{m,k,i/d}}.
\end{eqnarray*} \end{theorem}

\begin{proof} The left hand side is equal to

\begin{eqnarray*}
& & \sum_{w \in S_n} \sum_{m=0 \ mod \ n} Coef. \ of
\ q^m t^{k-cd(w)} \ in \ \frac{1}{(1-tq)\cdots(1-tq^n)} q^{maj(w)}
\prod x_i^{n_i(w)}\\
& = & \sum_{w \in S_n} \sum_{m=0 \ mod \ n} Coef. \ of \ q^m \ in \ \qb{n+k-cd(w)-1}{n-1}
q^{maj(w)} \prod x_i^{n_i(w)},
\end{eqnarray*} where the last step uses Theorem 349 on page 280 of \cite{HW}. Note by part 4 of Theorem \ref{typeAform} that this expression is precisely the cycle structure generating function under the measure $x_k$, multiplied by $k^{n-1}$.

	To complete the proof of the theorem, it must be shown that
the right hand side gives the cycle structure generating function for
degree $n$ polynomials over a field of $k$ elements with constant term
$1$. Let $\phi$ be a fixed generator of the multiplicative group of
the field $F_k$ of $k$ elements, and let $\tau_i$ be a generator of
the multiplicative group of the degree $i$ extension of $F_k$, with
the property that $\tau_i^{(q^i-1)/(q-1)}=\phi$. Recall Golomb's
correspondence between degree $i$ polynomials over $F_k$ and size $i$
aperiodic necklaces on the symbols $\{0,1,\cdots,k-1\}$. This
correspondence goes by taking any root of the polynomial, expressing
it as a power of $\tau_i$ and then writing this power base $k$ and
forming a necklace out of the coefficients of $1,k,k^2,\cdots,
k^{i-1}$. It is then easy to see that the norm of the corresponding
polynomial is $\phi$ raised to the sum of the necklace entries. The
result now follows from Lemma \ref{likeVic}. Note that there is no
$m=0$ term because the polynomial $z$ can not divide a polynomial with
constant term 1. \end{proof}

	Other evidence in favor of Conjecture 1 is its correctness for
type $C$ in odd characteristic. It is perhaps surprising that type $C$
is simpler to handle than type $A$. The reason for this phenomenon,
also encountered by Cellini \cite{Ce2}, is that for type $C$, the
coefficient of $w$ in the element $x_{2k+1}$ depends only on the
descent set of $w$ and not on the cyclic descent set.

	The following counting Lemma of \cite{FNP} will be helpful.

\begin{lemma} \label{countinvar} Let $e=1$ if $q$ is even and $e=2$ if
$q$ is odd. The number of monic, degree $n$ polynomials $f(z)$ over
$F_q$ with non-zero constant coefficient and invariant under the
involution $f(z) \mapsto f(0)^{-1} z^n f(\frac{1}{z})$ is

\[ \left\{ \begin{array}{ll}
																					e & \mbox{if
$n=1$}\\
																																			0 & \mbox{if $n$ \ is \ odd \ and \ $n>1$}\\
																						\frac{1}{n} \sum_{d|n \atop d \ odd} \mu(d) (q^{\frac{n}{2d}}+1-e)		& \mbox{Otherwise}
																																				\end{array}
			\right.			 \]

\end{lemma}
	
\begin{theorem} \label{checktypeC} Conjecture 1 holds for type $C$ in
odd characteristic. \end{theorem}

\begin{proof} From Lemma \ref{com}, when $q$ is odd the coefficient of
$w$ in $x_q$ is the same as the formula for $H_{C_n,q}(w)$ in
\cite{F3} (stated there for type $B$ but which is the same in type $C$
as it is defined only in terms of descents sets--for which types $B,C$
are equivalent--and not in terms of cyclic descents).

	The paper \cite{F3} showed, using delicate combinatorial
techniques of Reiner \cite{R}, that Conjecture 0 holds for type
$C$. Thus it suffices to show that the mass on a class $C$ of $W$
obtained by picking one of the $q^n$ semisimple conjugacy classes of
$Sp(2n,q)$ at random and applying $\Phi$ is equal to the mass on a
class $C$ of $W$ obtained by picking one of the $q^n$ semisimple
adjoint orbits of $Sp(2n,q)$ on its Lie algebra at random and then
applying $\Phi$.

	To do this, recall first that the semisimple conjugacy classes
of $Sp(2n,q)$ are monic degree $2n$ polynomials $f(z)$ with non-zero
constant term invariant under the involution sending $f(z)$ to
$\bar{f}(z)=\frac{z^{2n}f(\frac{1}{z})}{f(0)}$. The conjugacy classes
of $C_n$ correspond to pairs of vectors $(\vec{\lambda},\vec{\mu})$
where $\vec{\lambda}=(\lambda_1,\cdots,\lambda_n)$,
$\vec{\mu}=(\mu_1,\cdots,\mu_n)$ and $\lambda_i$ (resp. $\mu_i$) is
the number of positive (resp. negative) $i$ cycles of an element of
$C_n$, viewed as a signed permutation. From Section 3 of \cite{C2} and
Section 2 of \cite{SS}, the map $\Phi$ can be described as follows.
Factor $f$ uniquely into irreducibles as \[ \prod_{\{\phi_j,
\bar{\phi_j}\}} [\phi_j \bar{\phi_j}]^{r_{\phi_j}} \prod_{\phi_j :
\phi_j = \bar{\phi_j}} \phi_j^{s_{\phi_j}} \] where the $\phi_j$ are
monic irreducible polynomials and $s_{\phi_j} \in \{0,1\}$. Note that
all terms in the second product have even degree. This is because by
Lemma \ref{countinvar} all $\phi$ invariant under the involution
$\bar{}$ have even degree, except possibly $z \pm 1$. However $z \pm
1$ must each appear with even multiplicity as factors of the
characteristic polynomial of an element of $Sp(2n,q)$ and hence only
contribute to the first product. The class of $W$ corresponding to $f$
is then determined by setting $\lambda_i(f)=\sum_{\phi: deg(\phi)=i}
r_{\phi}$ and $\mu_i(f)=\sum_{\phi: deg(\phi)=2i} s_{\phi}$.

	Next, recall that the semisimple orbits of $Sp(2n,q)$ on its
Lie algebra are monic degree $2n$ polynomials $f(z)$ satisfying
$f(z)=f(-z)$. Arguing as in the previous paragraph, the description of
the map $\Phi$ is similar. One factors $f$ uniquely into irreducibles
as \[ \prod_{\{\phi_j(z), \phi_j(-z)\}} [\phi_j(z)
\phi_j(-z)]^{r_{\phi_j}} \prod_{\phi_j : \phi_j(z) = \phi_j(-z)}
\phi_j(z)^{s_{\phi_j}} \] where the $\phi_j$ are monic irreducible
polynomials and $s_{\phi_j} \in \{0,1\}$. Then
$\lambda_i(f)=\sum_{\phi: deg(\phi)=i} r_{\phi}$ and
$\mu_i(f)=\sum_{\phi: deg(\phi)=2i} s_{\phi}$. Thus the theorem will
follow if it can be shown that the number of monic, degree $2m$
irreducible polynomials satisfying $f=\bar{f}$ is equal to the number
of monic, degree $2m$ irreducible polynomials satisfying
$f(z)=f(-z)$. From Lemma \ref{countinvar}, the first quantity is \[
\frac{1}{2m} \sum_{d|m \atop d \ odd} \mu(d) (q^{\frac{m}{d}}-1).\]
This formula agrees with the second quantity, as computed in
\cite{R}. \end{proof}

\section{Dynamical Systems} \label{dynamical}

	This section is divided into two subsections. The first
indicates the relationship of Conjecture 1 for type $C$ even
characteristic to Gannon's enumeration of unimodal maps by cycle
structure. His results are then reexpressed in a form amenable to
asymptotic analysis and some conclusions are drawn. Along the way a
result from dynamical systems is given a combinatorial proof. The
second subsection considers enumeration of type $A$ Bayer-Diaconis
shuffles by the shape of their cycles (which is a finer invariant than
their length).

\subsection{Hyperoctahedral Shuffles} \label{subhyper}

	Theorem \ref{reduce} reformulates Conjecture 1 in type $C$
even characteristic.

\begin{theorem} \label{reduce} Let $\lambda_i(w)$ and $\mu_i(w)$ be
the number of positive and negative $i$-cycles of a permutation $w$ in
some $C_n$ and let $cd(w)$ be the number of cyclic descents of
$w$. Then in even characteristic Conjecture 1 follows from the
generating function identity

\[ \sum_{n \geq 0} \frac{u^n}{(1-t)^{n+1}} \sum_{w \in C_n} t^{cd(w)}
\prod_{i \geq 1} x_i^{\lambda_i(w)} y_i^{\mu_i(w)} = \sum_{s \geq 0}
t^s \prod_{m \geq 1} (\frac{1+x_mu^m}{1-y_mu^m})^{\frac{1}{2m}
\sum_{d|m \atop d \ odd} \mu(d) (2s)^{\frac{m}{d}}}.\] \end{theorem}

\begin{proof} Suppose that the identity of the theorem is true. Then
taking coefficients of $t^s$ on both sides would yield the equation
 
\[ 1+\sum_{n \geq 1} u^n \sum_{w \in C_n} {s+n-cd(w) \choose n} \prod_{i
\geq 1} x_i^{\lambda_i(w)} y_i^{\mu_i(w)} = \prod_{m \geq 1}
(\frac{1+x_m u^m}{1-y_m u^m})^{\frac{1}{2m} \sum_{d|m \atop d \ odd}
\mu(d) (2s)^{\frac{m}{d}}}.\] Set $s=\frac{q}{2}$ (which is an integer
since $q$ is assumed even). Taking the coefficient of $u^n \prod_i
x_i^{\lambda_i} y_i^{\mu_i}$ on the left hand side of this equation
and dividing by $q^n$ gives by Lemma \ref{com} the probability that an
$w$ chosen according to the $x_q$ probability measure is in a
conjugacy class with $\lambda_i$ positive $i$-cycles and $\mu_i$
negative $i$-cycles for each $i$. By Lemma \ref{countinvar}, doing the
same to the right hand side of the equation gives the probability that
when one factors a uniformly chosen random degree $2n$ polynomial over
$F_q$ as

\[ \prod_{\{\phi_j(z), \phi_j(-z)\}} [\phi_j(z)
\phi_j(-z)]^{r_{\phi_j}} \prod_{\phi_j : \phi_j(z) = \phi_j(-z)}
\phi_j(z)^{s_{\phi_j}} \] (with $\phi_j$ are monic irreducible
polynomials and $s_{\phi_j} \in \{0,1\}$), that one obtains
$\lambda_i=\sum_{\phi: deg(\phi)=i} r_{\phi}$ and $\mu_i=\sum_{\phi:
deg(\phi)=2i} s_{\phi}$. The theorem now follows from the fact that in
even characteristic the map $\Phi$ of Conjecture 1 has the same
description as in the proof of Theorem \ref{checktypeC}. \end{proof}

	We remark that a generating function identity very similar to
that in Theorem \ref{reduce} appears in \cite{R}. It must be the case
that a modification of Reiner's bijections will lead to a proof of the
identity.

	The next result relates to the enumeration of unimodal
permutations by cycle structure. Here a unimodal permutation $w$ on
the symbols $\{1,\cdots,n\}$ is defined by requiring that there is
some $i$ with $1 \leq i \leq n$ such that the following two properties
hold:
 \begin{enumerate} \item If $a<b\leq i$, then $w(a)<w(b)$.
  \item If $i \leq a<b$, then $w(a)>w(b)$.
  \end{enumerate} Thus $i$ is where the maximum is achieved, and the permutations
$12\cdots n$ and $nn-1\cdots 1$ are counted as unimodal. For each
fixed $i$ there are ${n-1 \choose i-1}$ unimodal permutations with
maximum $i$, hence a total of $2^{n-1}$ such permutations. Note that
this is exactly one half the number of terms in an inverse 2-shuffle
for type $C_n$.

	Motivated by biology and dynamical systems, Rogers \cite{Ro}
posed the problem of counting unimodal permutations by cycle
structure. This problem was solved by Gannon who gave a constructive
proof of the following elegant (and more fundamental) result. For its
statement, one defines the shape $s$ of a cycle $(i_1 \cdots i_k)$ on
some $k$ distinct symbols (call them $\{1,\cdots,k\}$) to be the cycle
$\{\tau(i_1) \cdots \tau(i_k)\}$ where $\tau$ is the unique order
preserving bijection between $\{i_1,\cdots,i_k\}$ and
$\{1,\cdots,k\}$.

\begin{theorem} \label{Gannon} (\cite{Ga}) Let $s_1,s_2,\cdots$ denote
the possible shapes of transitive unimodal permutations. Then the
number of unimodal permutations with $n_i$ cycles of shape $s_i$ is
$2^{l-1}$, where $l$ is the number of $i$ for which $n_i>0$. \end{theorem}

	Theorem \ref{Gannon} can be rewritten in terms of generating
functions.

\begin{cor} \label{cindex} Let $n_s(w)$ be the number of cycles of $w$
of shape $s$. Let $|s|$ be the number of elements in $s$. Then

\[ 1+\sum_{n=1}^{\infty} \frac{u^n}{2^{n-1}} \sum_{w \in S_n \atop w \
unimodal} \prod_{s \ shape} x_s^{n_s(w)} = \prod_{s \ shape}
(1+\frac{2x_su^{|s|}}{2^{|s|}-x_su^{|s|}})\]

\[ (1-u)+\sum_{n=1}^{\infty} \frac{(1-u) u^n}{2^{n-1}} \sum_{w \in S_n
\atop w \ unimodal} \prod_{s \ shape} x_s^{n_s(w)} = \prod_{s \ shape}
(\frac{2^{|s|}+x_su^{|s|}}{2^{|s|}+u^{|s|}})(\frac{2^{|s|}-u^{|s|}}{2^{|s|}
-x_su^{|s|}})\] \end{cor}

\begin{proof} For the first equation, consider the coefficient of
$\prod_s x_s^{n_s} u^{\sum |s| n_s}$ on the left hand side. It is the
probability that a uniformly chosen unimodal permutation on $\sum |s|
n_s$ symbols has $n_s$ cycles of shape $s$. The coefficient on the
right hand side is $2^{|\{s:n_s>0\}|-\sum n_s}$. These are equal by
Theorem \ref{Gannon}. To deduce the second equation, observe that
setting all $x_s=1$ in the first equation gives that \[ \frac{1}{1-u}
= \prod_{s \ shape} (1+\frac{2u^{|s|}}{2^{|s|}-u^{|s|}}) .\] Taking
reciprocals and multiplying by the first equation yields the second
equation. \end{proof}

{\bf Remarks:} \begin{enumerate}

\item The second equation in Corollary \ref{cindex} has an attractive
probabilistic interpretation. Fix $u$ such that $0<u<1$. Then choose a
random symmetric group so that the chance of getting $S_n$ is equal to
$(1-u)u^n$. Choose a unimodal $w \in S_n$ uniformly at random. Then
the random variables $n_s(w)$ are independent, each having
distribution a convolution of a binomial$(\frac{u^{|s|}}
{2^{|s|}+u^{|s|}})$ with a geometric$(1-\frac{u^{|s|}}{2^{|s|}})$.

	As another illustration of the second equation in Corollary
\ref{cindex}, we deduce the following corollary, extending the
asymptotic results in \cite{Ga} that asymptotically $2/3$ of all
unimodal permutations have fixed points and $2/5$ have $2$-cycles.

\begin{cor} \label{bigeom} In the $n \rightarrow \infty$ limit, the
random variables $n_s$ converge to the convolution of a
binomial$(\frac{1} {2^{|s|}+1})$ with a
geometric$(1-\frac{1}{2^{|s|}})$ and are asymptotically independent.
\end{cor}

\begin{proof} The result follows from the claim that if $f(u)$ has a
Taylor series around 0 and $f(1)<\infty$ and $f$ has a Taylor series
around 0, then the $n \rightarrow \infty$ limit of the coefficient of
$u^n$ in $\frac{f(u)}{1-u}$ is $f(1)$. To verify the claim, write the
Taylor expansion $f(u) = \sum_{n=0}^{\infty} a_n u^n$ and observe that
the coefficient of $u^n$ in $\frac{f(u)}{1-u} = \sum_{i=0}^n a_i$.
\end{proof}

	There is still more to be done, for instance studying the
length of the longest cycle.

\item The type $C_n$ shuffles also relate to dynamical systems in
another way, analogous to the type $A$ construction for Bayer-Diaconis
shuffles \cite{BD} described in the next subsection. Here we describe
the case $k=2$. One drops $n$ points in the interval $[-1,1]$
uniformly and indepedently. Then one applies the map $x \mapsto
2|x|-1$. The resulting permutation can be thought of as a signed
permutation, since some points preserve and some reverse
orientation. From Proposition \ref{rel}, this signed permutation
obtained after iterating this map $r$ times has the distribution of
the type $C_n$ shuffle with $k=2^r$. Lalley \cite{La1} studied the
cycle structure of random permutations obtained by tracking $n$
uniformly dropped points after iterating a map a large number of
times. His results applied to piecewise monotone maps, and he proved
that the limiting cycle structure is a convolution of
geometrics. Hence Corollary \ref{bigeom} shows that Lalley's results
do not extend to functions such as $x \mapsto 2|x|-1$.

\end{enumerate}

	Lemma \ref{countem} will be useful. Its proof by Rogers and
Weiss \cite{RW} used dynamical systems. Here a different proof is
presented using combinatorial machinery of Gessel and Reutenauer
\cite{G}.

\begin{lemma} \label{countem} The number of transitive unimodal
$n$-cycles is \[\frac{1}{2n} \sum_{d|n \atop d \ odd} \mu(d)
2^{\frac{n}{d}}.\] \end{lemma}

\begin{proof} Symmetric function notation from Chapter 1 of Macdonald
\cite{Mac} is used. Thus $p_{\lambda},h_{\lambda}$ $,
e_{\lambda},s_{\lambda}$ are the power sum, complete, elementary, and
Schur symmetric functions parameterized by a partition $\lambda$. From
Theorem 2.1 of \cite{G}, the number of $w$ in $S_n$ with a given cycle
structure and descent set $D$ is the inner product of a Lie character
$L_n=\frac{1}{n} \sum_{d|n} \mu(d) p_d^{\frac{n}{d}}$ and a Foulkes
character $F_{C(D)}$. From the proof of Corollary 2.4 of \cite{G},
$F_{C(D)}=\sum_{|\lambda|=n} \beta_{\lambda} s_{\lambda}$ where
$\beta_{\lambda}$ is the number of standard tableaux of shape
$\lambda$ with descent composition $C(D)$. Thus the sought number is
\[ <\frac{1}{n} \ \sum_{d|n} \mu(d) p_d^{\frac{n}{d}},
e_n+\sum_{i=2}^{n-1} s_{i,(1)^{n-i}}+h_n>.\] Expanding these Schur
functions using exercise 9 on page 47 of \cite{Mac}, using the fact
that the $p_{\lambda}$ are an orthogonal basis of the ring of
symmetric functions with known normalizing constants (page 64 of
\cite{Mac}), and using the expansions of $e_n$ and $h_n$ in terms of
the $p_{\lambda}$'s (page 25 of \cite{Mac}) it follows that

\begin{eqnarray*}
 & & <\frac{1}{n} \ \sum_{d|n} \mu(d) p_d^{\frac{n}{d}},
e_n+\sum_{i=2}^{n-1} s_{i,(1)^{n-i}}+h_n>\\
&=& <\frac{1}{n} \ \sum_{d|n} \mu(d) p_d^{\frac{n}{d}},
\sum_{i \ even} h_i e_{n-i}>\\
&=& \frac{1}{n} \sum_{d|n} \mu(d) <p_d^{\frac{n}{d}},\sum_{i=1,\cdots,\frac{n}{d} \atop di \ even} h_{di} e_{n-di}>\\
&=& \frac{1}{n} \sum_{d|n} \mu(d) <p_d^{\frac{n}{d}}, p_d^{\frac{n}{d}} \sum_{i=1,\cdots,\frac{n}{d} \atop di \ even} \frac{(-1)^{n-di-\frac{n}{d}+i}}{d^{\frac{n}{d}} i! (\frac{n}{d}-i)!}>\\
& = & \frac{1}{n} \sum_{d|n} \mu(d) (-1)^{n-\frac{n}{d}} \sum_{i=1,\cdots,\frac{n}{d} \atop di \ even} (-1)^i {\frac{n}{d} \choose i}\\
& = & \frac{1}{2n} \sum_{d|n \atop d \ odd} \mu(d) 2^{\frac{n}{d}}.
\end{eqnarray*}
\end{proof}

	Theorem \ref{re} relates Conjecture 1 to the enumeration of
unimodal permutations by cycle structure. Due to Corollary
\ref{cindex}, it can be taken as evidence in favor of Conjecture 1.

\begin{theorem} \label{re} Suppose that Conjecture 1 holds for type
$C$ when $k=2$. Let $n_i(w)$ be the number of $i$-cycles of $w \in
S_n$. Then

\[ 1+\sum_{n=1}^{\infty} \frac{u^n}{2^{n-1}} \sum_{w \in S_n \atop w \
unimodal} \prod_i x_i^{n_i(w)} = \prod_i
(\frac{2^i+x_iu^i}{2^i+u^i})(\frac{2^i-u^i}{2^i
-x_iu^i})^{\frac{1}{2i} \sum_{d|i \atop d \ odd} \mu(d)
2^{\frac{i}{d}}}.\] \end{theorem}

\begin{proof} From the first equation in the proof of Theorem
\ref{reduce} and Lemma \ref{countem}, it suffices to prove that

\[ 1+\sum_{n \geq 1} \frac{u^n}{2^n} \sum_{w \in C_n} {s+n-cd(w)
\choose n} \prod_{i \geq 1} x_i^{\lambda_i(w)+\mu_i(w)} =
1+\sum_{n=1}^{\infty} \frac{u^n}{2^{n-1}} \sum_{w \in S_n \atop w \
unimodal} \prod_i x_i^{n_i(w)} .\] For this it is enough to define a
$2$ to $1$ map $\eta$ from the $2^n$ type $C_n$ characteristic 2
shuffles to unimodal elements of $S_n$, such that $\eta$ preserves the
number of $i$-cycles for each $i$, disregarding signs.

	To define $\eta$, recalling Proposition \ref{rel} observe that
the $2$ shuffles are all ways of cutting a deck of size $n$, then
flipping the first half, and choosing a random interleaving. For
instance if one cuts a 12 card deck at position 6, such an
interleaving could be \[[-6,-5,7,8,-4,9,-3,10,-2,11,-1,12].\] Observe
that taking the inverse of this permutation and disregarding signs
gives \[[11,8,6,5,2,1,3,4,6,8,10,12].\] Next one conjugates by the
involution transposing each $i$ with $n+1-i$, thereby obtaining a
unimodal permutation. Note that this map preserves cycle structure,
and is $2$ to $1$ because the first symbol (in the example $-6$, can
always have its sign reversed yielding a possible shuffle).
\end{proof}

	We remark that Theorem \ref{re} gives a very natural
generalization of unimodal maps, namely the inverses of the elements
$x_k$, disregarding sign.

\subsection{Bayer-Diaconis Shuffles} \label{subBayer} 

	This subsection proves an analog of Theorem \ref{Gannon} for
the type $A$ riffle shuffles of \cite{BD}. First recall that their $k$-
shuffles on the symmetric group $S_n$ may be described as follows:

	Step 1: Start with a deck of $n$ cards face down. Choose $k$
numbers $j_1,\cdots,j_k$ multinomially with the probability of getting
$j_1,\cdots,j_k$ equal to $\frac{{n \choose
j_1,\cdots,j_k}}{k^n}$. Make $k$ stacks of cards of sizes
$j_1,\cdots,j_k$ respectively.

	Step 2: Drop cards from packets with probability proportional
to packet size at a given time. Equivalently, choose uniformly at
random one of the ${n \choose j_1,\cdots,j_k}$ interleavings of the
packets.

{\bf Remarks:}
\begin{enumerate}
\item Note that for any permutation obtained via a riffle shuffle,
there are at most $k$ rising sequences, corresponding to the pile
sizes. Thus these riffle shuffles can be thought of as a discrete
version of the map $x \mapsto kx$ mod $1$, which has exactly $k$
equally sized intervals on which it is monotonically increasing.

\item The following observation from \cite{BD} is useful and shows
that their type $A$ shuffle has a rigorous interpretation in terms of
dynamical systems. Namely if one drops $n$ points in the interval
$[0,1]$ uniformly and independently and then applies the map $x
\rightarrow kx \ mod \ 1$, the (random) permutation describing the
reordering of the points is exactly a $k$ shuffle on $S_n$.
\end{enumerate}

	The type $A$ shuffles of this subsection were studied by cycle
length in the paper \cite{DMP}. Motivated by Theorem \ref{Gannon}, we
consider their distribution by shape. The key tool is a bijection of
Gessel and Reutenauer \cite{G} which is analogous to Gannon's monoidal
construction \cite{Ga}. Let us review this bijection.

 	Define a necklace on an alphabet to be a sequence of
cyclically arranged letters of the alphabet. A necklace is said to be
primitive if it is not equal to any of its non-trivial cyclic
shifts. For example, the necklace $(a\ a\ b\ b)$ is primitive, but the
necklace $(a\ b\ a\ b)$ is not. Given a word $w$ of length $n$ on an
ordered alphabet, the 2-row form of the standard permutation $st(w)
\in S_n$ is defined as follows. Write $w$ under $1 \cdots n$ and then
write under each letter of $w$ its lexicographic order in $w$, where
if two letters of $w$ are the same, the one to the left is considered
smaller. For example (page 195 of $\cite{G}$):

\[ \begin{array}{c c c c c c c c c c c c c c}
	& & 1 & 2 & 3 & 4 & 5 & 6 & 7 & 8 & 9 & 10 & 11 & 12\\
	w & = & b & b & a & a & b & c & c & c & b & c & b & b\\
	st(w) &= & 3 & 4 & 1 & 2 & 5 & 9 & 10 & 11 & 6 & 12 & 7 & 8
	\end{array} \]

	For a finite ordered alphabet $A$, Gessel and Reutenauer
(loc. cit.) give a bijection $U$ from the set of length $n$ words $w$
of onto the set of finite multisets of necklaces of total size $n$,
such that the cycle structure of $st(w)$ is equal to the cycle
structure of $U(w)$. To define $U(w)$, one replaces each number in the
necklace of $st(w)$ by the letter above it. In the example, the
necklace of $st(w)$ is $(1\ 3), (2\ 4), (5), (6\ 9), (7\ 11\ 8\ 12\
10)$. This gives the following multiset of necklaces on $A$:

\[ (a\ b)(a\ b)(b)(b\ c)(b\ c\ b\ c\ c) \] 

\begin{prop} \label{analog} Let $s_1,s_2,\cdots$ denote the possible
shapes of cycles of $k$-shuffles. Let $d(s^{-1})$ be the number of
descents of the inverse of the cycle of shape $s$. Then the
probability that a $k$-shuffle has $n_s$ cycles of shape $s$ is \[
\frac{1}{k^n} \prod_{s \ shape} {{|s|+k-d(s^{-1})-1 \choose |s|} + n_s
-1 \choose n_s}.\] \end{prop}

\begin{proof} From the above bijection, the $k^n$ inverse shuffles
correspond to multisets of primitive necklaces on the symbols
$\{1,\cdots,k\}$ in a cycle length preserving way. Thus the sought
probability is $\frac{1}{k^n}$ times the product (over shapes $s$) of
the number of ways choosing a multiset of size $n_s$ from the
primitive necklaces on $|s|$ symbols and of shape $s$. The number of
such necklaces is ${|s|+k-d(s)-1 \choose |s|}$, because from \cite{BD},
this is the probability that an inverse $k$ shuffle on $|s|$ symbols
gives the cycle $s$. \end{proof}

	Proposition \ref{analog} can be rewritten in terms of
generating functions as

\[ 1+\sum_{n=1}^{\infty} \frac{u^n}{k^n} \sum_{w \in S_n} {n+k-d(w)-1
\choose n} \prod_{s \ shape} x_s^{n_s(w)} = \prod_{s \ shape}
(\frac{1}{1-u^{|s|}x_s})^{{|s|+k-d(s^{-1})-1 \choose |s|}}.\]

\section{Refining the Map $\Phi$ to the Weyl Group} \label{refine}

	This section discusses the possibility of refining the map
$\Phi$ so as to naturally associate to a semisimple conjugacy class
$c$ an element $w$ of the Weyl group. Furthermore the conjugacy class
of $w$ should be $\Phi(c)$, and the induced probability measure on $W$
should agree with that given by $x_q$. Such a result could be
important in algebraic number theory, as semisimple conjugacy classes
play a key role in the Langlands program. For example consider a
simple algebraic extension of $Q$ whose generator has minimal
polynomial $f(x)$. At unramified primes $p$, the conjugacy class of
the Frobenius automorphism, viewed as a permutation of the roots, is
given by the degrees of the irreducible factors of the mod $p$
reduction of $f$. This is exactly the map $\Phi$ in type $A$. It is
standard to encode such data over all primes into a generating
function (see \cite{Gel} for a survey). The hope is that a refinement
of $\Phi$ will lead to refined number theoretic constructs.

	This program of refining $\Phi$ was partly successful in the
case of semisimple orbits on the Lie algebra \cite{F3}. There the key
idea was a combinatorial bijection of Gessel, combined with elementary
Galois theory. Those ideas apply to the $C_n$ conjugacy class setting
in this paper, in odd characteristic. Here, however, we wish to pursue
another possibility, using the geometry of the affine Weyl group. A
good background reference is Section 3.8 of Carter \cite{C}. For the
next paragraph we follow his treatment.

	Let $Y$ be the coroot lattice and $W$ the Weyl group, so that
$<Y,W>$ is the affine Weyl group. The group $<Y,W>$ acts on the vector
space $Y \otimes R$ with $Y$ acting by translations $v \mapsto v+y$
and $W$ acting by orthogonal transformations. The affine Weyl group
has a fundamental region in $Y \otimes R$ given by \[ \bar{A} = \{v
\in Y \otimes R| <\alpha_i,v> \geq 0 \ for \ i=1,\cdots,r,
<\alpha_0,v> \leq 1\}.\] Each element of $Y \otimes R$ is equivalent
to exactly one element of the fundamental chamber. Let $Q_{p'}$ be the
additive group of rational numbers $\frac{s}{t}$ where $s,t \in Z$ and
$t$ is not divisible by $p$ (the characteristic). Proposition 3.8.1 of
\cite{C} shows that there is an action of $F$ on $\bar{A_{p'}}=\bar{A}
\cap (Y \otimes Q_{p'})$ given by taking the image of $v \in
\bar{A_{p'}}$ to be the unique element of $\bar{A}$ equivalent to
$F(a)$ under $<Y,W>$. There are $q^r$ elements of $\bar{A_{p'}}$
stable under this action, and by Proposition 3.7.3 of \cite{C}, they
correspond to the semisimple conjugacy classes of $G^F$.

	Thus each stable point $v \in \bar{A_{p'}}$ satisfies
$F(v)=wv+y$ for some $w \in W, y \in Y$. Considering the set of $w \in
W$ for which there is a $y \in Y$ satisfying this equation, the idea
is to choose such a $w$ uniquely. The most natural possibility is to
pick $w$ of minimal length. It is not clear in general that such a $w$
is unique. Lemma \ref{unique} shows that such a $w$ is unique for type
$A$, and further that it has the same conjugacy class as that obtained
by the map $\Phi$. The preparatory Lemma \ref{clas} will be helpful.

\begin{lemma} \label{clas} Let $w'$ be an element of the symmetric
group $S_n$. Let $I_1,\cdots,I_t$ be such that $\{1,\cdots,n\}$ is
their disjoint union. Suppose that $w'$ permutes $I_1,\cdots,I_t$, and
furthermore does so in such a way that if $j \in I_k$ is in an
$i$-cycle under $w'$, then $I_k$ is in an $i$-cycle in the action of
$w'$ on the set $\{I_1,\cdots,I_t\}$. Then $w'$ is conjugate to the
unique minimal length coset representative of the coset $w' (S_{I_1}
\times \cdots \times S_{I_t})$.\end{lemma}

\begin{proof} Let $w$ be the unique minimal length coset
representative of the coset $w'(S_{I_1} \times \cdots \times
S_{I_t})$. Then $w$ is given explicitly by reordering the images under
$w'$ of the elements within each $I_1,\cdots,I_t$ to be in increasing
order. Consider $j$ in some $I_k$. Suppose that $j$ is the $rth$
smallest element in $I_k$. If the $w'$-orbit of $I_k$ is
$I_{l_1},\cdots,I_{l_i}$, then the $w$-orbit of $j$ consists of the
$r$th smallest elements of each of $I_{l_1},\cdots,I_{l_i}$, and hence
has size $i$. By hypothesis, the $w'$ orbit of $j$ also has size $i$,
implying the result.  \end{proof}

\begin{lemma} \label{unique} In type $A$, for each stable point $v \in
\bar{A_{p'}}$, there is a unique shortest $w$ such that $F(v)=wv+y$
for some $y \in Y$. Letting $c$ be the semisimple conjugacy class of
$G^F$ corresponding to $v$, the conjugacy class of $w$ is
$\Phi(c)$.\end{lemma}

\begin{proof} Denoting the coordinates of $v$ as $x_1,\cdots,x_n$,
recall that $x_1 \geq x_2 \geq \cdots \geq x_n$, $x_1-x_n \leq 1$, and
that $x_1+\cdots+x_n=0$. Semisimple conjugacy classes $c$ of $SL(n,q)$
are monic polynomials with constant term 1. Since the roots of an
irreducible degree $i$ polynomial $\phi$ are a Frobenius orbit of some
element of a degree $i$ extension of $F_q$ lying in no smaller
extension, these roots correspond to some subset of values of
$\{x_1,\cdots,x_n\}$ (there may be repetition among the $x$'s) which
are of the form $\{\frac{a_{j_1}}{q^i-1}, \cdots,
\frac{a_{j_i}}{q^i-1}\}$ with $a_{j_1},\cdots,a_{j_i}$ integers. For
clarity of exposition, we point out that the roots of $\phi$ are
simply $\tau^{a_{j_1}},\cdots,\tau^{a_{j_i}}$ for $\tau$ is a
generator of the multiplicative group of the degree $i$ extension of
$F_q$ and was determined by the bijection in Proposition 3.7.3 of
\cite{C}.

	The $F$-action performs a cyclic permutation of these
$i$-values. Letting $r_i$ be such that the class $c$ corresponds to a
polynomial whose factorization into irreducibles has $r_i$ factors of
degree $i$, it follows that there is a permutation $w'$, with $r_i$
$i$-cycles, which satisfies the equation $F(v)=w'v+y$ for some $y \in
Y$ and such that $\Phi(c)$ is the conjugacy class of $w'$. Now suppose
that some other $w$ satisfies the equation $F(v)=wv+z$ for some $z \in
Y$. Then $w'v-wv=y-z$, which implies that $w'v=wv$ where equality
means equality of coordinates mod 1. Let $I_1,\cdots,I_t$ be the
partition of $\{1,\cdots,n\}$ induced by the equivalence relation that
$i \sim j$ when $x_i=x_j$ mod 1. Then $(w')^{-1}w$ is contained in the
parabolic subgroup $S_{I_1} \times \cdots S_{I_t}$. Thus the set of
all permutations $w$ satisfying the equation $F(v)=wv+y$ for some $y$
in $Y$ is simply the coset $w' (S_{I_1} \times \cdots S_{I_t})$. Lemma
\ref{clas} implies that the unique minimal length coset representative
is conjugate to $w'$. \end{proof}

	Theorem \ref{spec} verifies that the suggested construction
succeeds for the symmetric group $S_3$.

\begin{theorem} \label{spec} In type $A_2$, pick one of the $q^2$
$F$-stable elements $v$ of $\bar{A_{p'}}$ uniformly at
random. Let $w$ be the shortest element of $W$ such that the equation
$q v =w v+y$ has a solution for some $y \in
Y$. The induced probability measure on $W$ agrees with that given from
$x_q$.  \end{theorem}

\begin{proof} First it will be shown that the probability of the
identity element agrees for both measures. Fixed points under the
action of $F$ (which is multiplication by $q$) correspond to solutions
to the equations $\frac{c_1}{q-1}+\frac{c_2}{q-1}+\frac{c_3}{q-1}=0$,
$\frac{c_1}{q-1} \geq \frac{c_2}{q-1} \geq \frac{c_3}{q-1}$,
$\frac{c_1}{q-1}-\frac{c_2}{q-1} \leq 1$. The number of such solutions
is \[ \sum_{I \subseteq \tilde{\Pi}} a_{q-1,I}, \] so they can be
counted by the methods of Theorem \ref{checkS3}. Doing this, and
comparing with the formula for $x_q$ in Theorem \ref{checkS3} proves
the result for the identity element.

	Next we argue that the measures agree for the permutation
whose 2-line form is $321$ (i.e. the longest element). From Theorem
\ref{checkS3}, the mass which $x_q$ puts on this element is seen to be

\[ \begin{array}{ll} \frac{q^2-q}{6q^2} & \mbox{if $q=0,1 \ mod \ 3$}\\
\frac{q^2-q+4}{6q^2} & \mbox{if $q=2 \ mod \ 3$} \end{array} \] Now
consider an $F$-stable element of $\bar{A_{p'}}$ which is stabilized
by $321$ but by no other elements. Next we argue that necessary and
sufficient conditions for the coordinates $(c_1,c_2,c_3)$ of such a
point are the inequalities \[ c_1=\frac{aq+b}{q^2-1} \geq
c_2=\frac{(q^2-1)-(a+b)(q+1)}{q^2-1} \geq
c_3=\frac{bq+a-(q^2-1)}{q^2-1} \] and $c_1-c_3 \leq 1$ with $a,b$
integers satisfying $(q-1) \geq b > a \geq 0$. Here $c_2=-(c_1+c_3)$,
so it is enough to argue for the form of $c_1$ and $c_3$. Observe that
$0<c_1<1$ and $c_3<0$, since the origin, though stabilized by $321$,
is also stabilized by the identity, which is shorter. The form of
$c_1$ follows since $F^2$ stabilizes $(c_1,c_2,c_3)$, but $F$ does not
(note that $b \geq a$ since $c_1-c_3 \leq 1$ and $b \neq a$ since $F$
does not stabilize $(c_1,c_2,c_3)$).

	The inequalities $c_1 \geq c_2 \geq c_3$ and $c_1 - c_3 \leq
1$ imply without much difficulty that $2a+b \geq q-1$. Rewriting $c_1
\geq c_2$ and $c_2 \geq c_3$ gives that

\[ q(2a+b-(q-1))+(a+2b-(q-1)) \geq 0 \] \[ q^2-1 \geq q(2b+a-(q-1) +
((2a+b)-(q-1)) \] Viewing the sums in these inequalities as base $q$
expansions and using the facts that $2a+b \geq q-1$ and $2b+a \geq
2a+b$, one concludes that our system of inequalities is equivalent to
the system of 3 inequalities

\[q-1 \geq b > a \geq 0\] \[ 0 \leq 2b+a \leq 2(q-1)\] \[ 2a+b \geq
q-1\] To count the number of integer solutions, we first count
solutions ignoring the third constraint, then subtract off solutions
satisfying $2a+b<q-1$. We only work out the case $q=2$ mod 3, the
other cases being similar. The number of solutions ignoring the third
inequality is 

\[ \sum_{b=1}^{\frac{2q-1}{3}} b + \sum_{b=\frac{2q-1}{3}+1}^{q-1}
(2q-2b-1) = \frac{2q^2-2q+2}{6}.\] The number of solutions obtained by
replacing the third inequality with $2a+b<q-1$ is readily computed to
be $\frac{q^2-q-2}{6}$, implying that the number of solutions to the
original inequalities is $\frac{q^2-q+4}{6}$, as desired.

	From Theorem \ref{checkS3}, the element $x_q$ assigns equal
mass to the remaining two tranpositions $213$ and $132$ in
$S_3$. Noting that $(c_1,c_2,c_3) \in \bar{A_{p'}}$ is $F$-stable if
and only if $(-c_3,-c_2,-c_1)$ is $F$-stable, it follows that the
construction of this section also assigns equal mass to $213$ and
$132$. Thus it remains to show that the construction of this section
and $x_q$ assign equal mass to the conjugacy class of all
transpositions. From Lemma \ref{unique}, the construction of this
section assigns equal mass to the conjugacy class of transpositions as
does the map $\Phi$. The result now follows by Theorem \ref{checkS3}.

	To argue for the class of 3-cycles, from Theorem \ref{checkS3}
the element $x_q$ assigns equal mass to the permutations $(123)$ and
$(132)$. Noting that $(c_1,c_2,c_3) \in \bar{A_{p'}}$ is $F$-stable if
and only if $(-c_3,-c_2,-c_1)$ is $F$-stable, it follows that the
construction of this section also assigns equal mass to $(123)$ and
$(132)$. Now argue as in the preceeding paragraph.
\end{proof}

	Proposition \ref{spec2} shows that the idea behind Theorem
\ref{spec} carries over to the identity conjugacy class in type $A$.

\begin{prop} \label{spec2} In type $A_{n-1}$, pick one of the
$q^{n-1}$ $F$-stable elements $v$ of $\bar{A_{p'}}$ uniformly at
random. Let $w$ be the shortest element of $W$ such that the equation
$q v =w v+y$ has a solution for some $y \in Y$. Then the chance that
$w$ is the identity is equal to the coefficient of the identity in
$x_q$. \end{prop}

\begin{proof} The proposition follows from Lemma \ref{unique} and
Corollary \ref{checkident}. \end{proof}

	Based on Lemma \ref{unique}, Theorem \ref{spec}, and Proposition
\ref{spec2}, we propose

{\bf Conjecture 2}: In type $A_n$, pick one of the $q^n$ $F$-stable
elements $v$ of $\bar{A_{p'}}$ uniformly at random. Let $w$ be the
shortest element of $W$ such that the equation $q v =w v+y$ has a
solution for some $y \in Y$. The induced probability measure on $W$
agrees with that given from $x_q$.

{\bf Problem:} Is the obvious analog of Conjecture 2 true in all
types?

\section{Acknowledgements} This research was supported by an NSF
Postdoctoral Fellowship.

\end{document}